\input amstex
 \documentstyle{amsppt}
 \magnification 1200
 \NoBlackBoxes
 \pagewidth{6.5 true in} 
 \pageheight{9.25 true in}

 \topmatter
 \title 
Lower bounds for moments of $L$-functions:  Symplectic and Orthogonal examples
 \endtitle
 \author 
 Z. Rudnick and K. Soundararajan
 \endauthor 
 \rightheadtext{Lower bounds for moments}
 \address School of Mathematical Sciences, Tel Aviv University,
Tel Aviv 69978, Israel
\endaddress
\email{rudnick{\@}post.tau.ac.il}
\endemail
\address{Department of Mathematics, University of Michigan, Ann Arbor,
Michigan 48109, USA} \endaddress \email{ksound{\@}umich.edu}
\endemail

\thanks{The first author is partially supported by a grant from the
  Israel Science Foundation. 
%Grant No. 2002088 from the
%United States-Israel Binational Science Foundation (BSF),
%Jerusalem, Israel. 
The second author is partially supported by the
National Science Foundation.}
\endthanks

 \endtopmatter
 \def\sumf{\sideset\and^{\flat} \to \sum}
 \def\sumstar{\sideset \and^{*} \to \sum}
 \def\lam{\lambda}
 \def\Lam{\Lambda}
 \def\Sym{\operatorname{Sym}}
 \def\sumh{\sideset \and^{h} \to \sum} 
\document 

 \head 1. Introduction \endhead
 
 \noindent An important problem in number theory asks for asmptotic 
 formulas for the moments of central values of $L$-functions varying in 
 a family.  This problem has been intensively studied in recent years, and 
 thanks to the pioneering work of  Keating and  Snaith [7], and the subsequent 
 contributions of  Conrey, Farmer, Keating, Rubinstein and Snaith [1], 
and Diaconu, Goldfeld and  Hoffstein [3] there are now well-established
conjectures for these moments.   The conjectured asymptotic formulas 
take different shapes depending on 
the symmetry group attached to the family of $L$-functions, given in the work of 
Katz and Sarnak  [6],  with 
three classes of formulas depending on whether the group in question is
unitary, orthogonal or symplectic.   While there are many known examples 
of asymptotic formulas dealing with the first few moments of a family 
of $L$-functions, in general the moment conjectures seem formidable.  
In [8] we recently gave a simple method to obtain lower bounds of 
the conjectured order of magnitude in many families of $L$-functions.  
In [8] we illustrated our method by working out lower bounds for 
$\sumstar_{\chi \pmod q} |L(\frac 12,\chi)|^{2k}$ where $q$ is a large
 prime, and the  
sum is over the primitive Dirichlet $L$-functions $\pmod q$.   This was an 
example of a `unitary' family of $L$-functions, and in this paper we round out 
the picture by providing lower bounds for moments of $L$-functions arising 
from orthogonal and symplectic families.

As our first example, we consider ${\Cal H}_k$ the set 
of Hecke eigencuspforms of weight $k$ for the 
full modular group $SL(2,{\Bbb Z})$.   We will 
think of the weight $k$ as being large, and note that 
${\Cal H}_k$ contains about $k/12$ forms.  Given $f\in {\Cal H}_k$ we 
write its Fourier expansion as 
$$
f(z) = \sum_{n=1}^{\infty} \lam_f(n) n^{\frac {k-1}{2}} e(nz), 
$$
where we have normalized the  Fourier coefficients so that 
the Hecke eigenvalues $\lam_f(n)$ satisfy Deligne's 
bound $|\lam_f(n)| \le \tau(n)$ where 
$\tau(n)$ is the number of divisors of $n$.  
Consider the associated $L$-function
$$
L(s,f) = \sum_{n=1}^\infty \lambda_f(n) n^{-s} =
\prod_p (1-\lambda_f(p) p^{-s} +p^{-2s})^{-1},  
$$
which converges absolutely in Re $(s)>1$ and extends 
analytically to the entire complex plane.  Recall that 
$L(s,f)$ satisfies the functional equation
$$
 \Lam(s,f):=(2\pi )^{-s}\Gamma(s+\tfrac{k-1}{2})L(s,f) =
i^k \Lam(1-s, f).   
$$
If $k\equiv 2\pmod 4$ then the sign of the functional equation is negative and
so $L(\tfrac 12,f)=0$.  We will therefore assume that $k\equiv 0 \pmod 4$.

While dealing with moments of $L$-functions in ${\Cal H}_k$, 
it is convenient to use the natural `harmonic weights' that 
arise from the Petersson norm of $f$.  Define the weight 
$$
\omega_f : =  \frac{(4\pi)^{k-1}}{\Gamma(k-1)} \langle f,f\rangle =
\frac{k-1}{2\pi^2} L(1,\Sym^2 f), 
$$
where $\langle f,f\rangle$ denotes the Petersson inner product.  
For a typical $f$ in ${\Cal H}_k$ the harmonic weight $\omega_f$ 
is of size about $k/12$, and so $\sum_{f\in {\Cal H}_k} \omega_f^{-1}$
is very nearly $1$.  The weights $\omega_f$ arise naturally in 
connection with the Petersson formula, and the facts mentioned above are 
standard and may be found in Iwaniec [4].  

For a positive integer $r$, we are interested in 
the $r$-th moment 
$$
\sumh_{f\in {\Cal H}_k } L(\tfrac 12, f)^r := 
\sum\Sb f\in {\Cal H}_k \endSb \frac{1}{\omega_f} L(\tfrac 12, f)^r. 
$$ 
This family of $L$-functions is expected to be of `orthogonal type' and 
the Keating-Snaith conjectures predict that for any given $r \in {\Bbb N}$ 
as $k\to \infty$ with $k\equiv 0\pmod 4$ we have
$$
\sumh_{f\in {\Cal H}_k } L(\tfrac 12, f)^r 
\sim C(r) (\log k)^{r(r-1)/2}, 
$$
for some positive constant $C(r)$.  This conjecture can be verified 
for $r=1$ and $r=2$, and if we permit an additional averaging 
over the weight $k$ then for $r=3$ and $4$ also.

\proclaim{Theorem 1}  For any given even natural number $r$, 
and weight $k\ge 12$ with $k\equiv 0 \pmod 4$, we have 
$$
\sumh_{f\in {\Cal H}_k} L(\tfrac 12, f)^r \gg_r (\log k)^{r(r-1)/2}.
$$
 \endproclaim

In fact, with more effort our method could be adapted to 
give lower bounds as in Theorem 1 for all rational numbers 
$r \ge 1$, rather than just even integers. 

Our other example involves the family of quadratic Dirichlet $L$-functions.
Let $d$ denote a fundamental discriminant, and let $\chi_d$ 
denote the corresponding real primitive character 
with conductor $|d|$.  We are interested in the class 
of quadratic Dirichlet $L$-functions $L(s,\chi_d)$.  Recall 
that, with ${\frak a}=0$ or $1$ depending on whether $d$ is 
positive or negative, these $L$-functions satisfy the 
functional equation 
$$
\Lam(s,\chi_d) := \Big(\frac{q}{\pi}\Big)^{\frac{s+{\frak a}}{2}} 
\Gamma\Big(\frac{s+{\frak a}}{2} \Big) L(s,\chi_d) 
= \Lam(1-s,\chi_d). 
$$
Notice that the sign of the functional equation is always 
positive, and it is expected that the central values 
$L(\tfrac 12,\chi_d)$ are all positive although this remains 
unknown.   
%For simplicity, we restrict the discussion to fundamental discriminants of the
%form $8d$ where $d$ is an odd, positive, square-free integer so that
%$\chi_{8d}(n)=(\frac{8d}n)$ is an even, real, primitive  character of
%conductor $8d$. We denote by $\sumf_{d\sim X}$ the sum over all such
%discriminants with $X/2<8d<X$. The number of such discriminants is
%asymptotically a constant multiple of $X$. 
This family is expected to be of  `symplectic' type and 
the Keating-Snaith conjectures predict that for any 
given $k \in {\Bbb N}$ and as $X\to \infty$  we have 
$$
\sumf_{|d|\le X}  L(\tfrac 12,\chi_d)^k \sim D(k) X (\log X)^{k(k+1)/2},
$$
for some positive constant $D(k)$, 
where the $\flat$ indicates that the sum is over fundamental discriminants.  
Jutila [5] established asymptotics for the first two moments of this 
family, and the third moment was evaluated in Soundararajan [10].

\proclaim{Theorem 2}  For every even natural number $k$ 
we have 
$$
\sumf_{|d| \le X} L(\tfrac 12,\chi_d)^k \gg_k X (\log X)^{k(k+1)/2}.
$$
\endproclaim 

As with Theorem 1, our method can be used to obtain lower bounds for these moments 
for all rational numbers $k$, taking care to replace $L(\tfrac 12,\chi_d)^k$ by $|L(\tfrac 12,\chi_d)|^k$ 
when $k$ is not an even integer.  In the case of the fourth moment we are able 
to get a lower bound $\ge (D(4)+o(1)) X (\log X)^{10}$, which matches exactly  
the asymptotic conjectured by Keating and Snaith.  The details of this calculation 
will appear elsewhere.

\head 2.  Proof of Theorem 1 \endhead

\noindent Let $x:= k^{\frac{1}{2r}}$ and consider 
$$
A(f): =A(f,x) = \sum_{n\le x} \frac{\lam_f(n)}{\sqrt{n}}.
$$
We will consider 
$$
S_1 := \sumh_{f\in {\Cal H}_k} L(\tfrac 12,f) A(f)^{r-1}, 
\ \ \text{ and } \ \ S_2: = 
\sumh_{f\in {\Cal H}_k} A(f)^r.
$$
Then H{\" o}lder's inequality gives, keeping 
in mind that $r$ is even so that $|A(f)|^r=A(f)^r$, 
$$
\Big(\sumh_{f\in {\Cal H}_k} L(\tfrac 12,f) A(f)^{r-1} \Big)^r 
\le \Big( \sumh_{f\in {\Cal H}_k } L(\tfrac 12,f)^r \Big) \Big( \sumh_{f\in {\Cal H}_k} 
A(f)^r \Big)^{r-1},
$$
so that 
$$
\sumh_{f\in {\Cal H}_k} L(\tfrac 12,f)^r \ge \frac{S_1^r}{S_2^{r-1}}. 
$$
We will prove Theorem 1 by finding the asymptotic orders of 
magnitude of $S_1$ and $S_2$.  

We begin with $S_2$.  To evaluate this, we must expand out 
$A(f)^r$ and group terms using the Hecke relations.  
To do this conveniently, let us denote by ${\Cal H}$ 
the ring generated over the integers by symbols $x(n)$ ($n\in {\Bbb N}$)
subject to the Hecke relations 
$$
x(1) =1, \qquad \text{and}\qquad x(m)x(n) = 
\sum_{d|(m,n)} x\Big(\frac{mn}{d^2}\Big).
$$
Thus ${\Cal H}$ is a polynomial ring on $x(p)$ where 
$p$ runs over all primes.  Using the Hecke relations we 
may write 
$$
x(n_1)\cdots x(n_r) = \sum_{t| \prod_{j=1}^{r} n_j } b_t(n_1,\ldots,n_r) 
x(t), 
$$
for certain integers $b_t(n_1,\ldots,n_r)$.  Note that $b_t(n_1,
\ldots,n_r)$ is symmetric in the variables $n_1$, $\ldots$, $n_r$, 
and that $b_t(n_1,\ldots,n_r)$ is always non-negative, and 
finally that $b_t(n_1,\ldots,n_r) \le \tau(n_1)\cdots \tau(n_r) 
\ll (n_1 \cdots n_r)^{\epsilon}$.  Of special importance for 
us will be the coefficient of $x(1)$ namely $b_1(n_1,\ldots,n_r)$.  
It is easy to see that $b_1$ satisfies a multiplicative 
property: if $(\prod_{j=1}^r m_j , \prod_{j=1}^{r} n_j) =1$ then 
$$
b_1(m_1n_1,\ldots, m_r n_r) = b_1(m_1,\ldots,m_r) b_1(n_1,\ldots,n_r).
$$
Thus it suffices to understand $b_1$ when the $n_1$, $\ldots$, $n_r$ 
are all powers of some prime $p$.  Here we note that $b_1(p^{a_1},
\ldots,p^{a_r})$ is independent of $p$, always lies between $0$ and
$(1+a_1)\cdots(1+a_r)$, and that it equals $0$ if $a_1+\ldots+a_r$ is 
odd.  If we write 
$$
B_r(n) = \sum\Sb n_1, \ldots, n_r \\ n_1 \cdots n_r =n \endSb 
b_1(n_1,\ldots,n_r), 
$$
then we find that $B_r(n)$ is a multiplicative function, that 
$B_r(n)=0$ unless $n$ is a square, and that $B_r(p^a)$ is 
independent of $p$ and grows at most polynomially in $a$.  Finally, and 
crucially, we note that 
$$
B_r(p^2) = r(r-1)/2,
$$
which follows upon noting that $b_1(p^2,1,\ldots,1)=0$ and 
that $b_1(p,p,1,\ldots,1)=1$.  

Returning to $S_2$ note that 
$$
A(f)^r = \sum_{n_1,\ldots,n_r \le x} 
\frac{1}{\sqrt{n_1\cdots n_r}} \sum_{t|n_1 \cdots n_r} b_t(n_1,\ldots,n_r)
\lam_f(t), 
$$
and so we require knowledge of $\sumh_{f\in {\Cal H}_k} \lam_f(t)$.  
This follows easily from Petersson's formula.  

\proclaim{Lemma 2.1} If $k$ is large, and $t$ and $u$ 
are natural numbers with 
$tu\le k^2/10^4$ then 
$$
\sumh_{f\in {\Cal H}_k} \lam_f(t)\lam_f(u) 
= \delta(t,u) + O(e^{-k}),
$$
where $\delta(t,u)$ is $1$ if $t=u$ and is $0$ otherwise. 
\endproclaim  
\demo{Proof} Petersson's formula (see [4]) gives
$$
\sumh_{f\in {\Cal H}_k} \lam_f(t)\lam_f(u) = \delta(t,u) + 2\pi i^k 
\sum_{c=1}^{\infty} \frac{S(t,u;c)}{c} J_{k-1} \Big(\frac{4\pi\sqrt{tu}}{c}
\Big). 
$$
Note that if $z\le 2k$ then $(z/2)^{k-1+\ell}/\Gamma(k-1+\ell) \le (z/2)^{k-1}/\Gamma(k-1)$ 
for all $\ell \ge 0$.  We now use the series representation for $J_{k-1}(z)$ 
which gives, for $z\le 2k$, 
$$
|J_{k-1}(z)| = \Big|\sum\Sb \ell=0\endSb^{\infty} 
\frac{(-1)^\ell (z/2)^{\ell}}{\ell !} \frac{(z/2)^{\ell+k-1}}{\Gamma(\ell+k-1)} 
\Big| \le \frac{(z/2)^{k-1}}{\Gamma(k-1)} e^{z/2}.
$$
Therefore, for $tu\le k^2/10^4$, we deduce that 
$$
\Big| J_{k-1}\Big(\frac{4\pi \sqrt{tu}}{c}\Big)\Big| 
\le \frac{(2\pi k/(100c))^{k-1}}{(k-2)!} e^{\pi k/50}.
$$
Using the trivial bound $|S(t,u;c)| \le c$ we 
conclude that 
$$
2\pi i^k \sum_{c=1}^{\infty} \frac{S(t,u;c)}{c} 
J_{k-1}\Big(\frac{4\pi \sqrt{tu}}{c}\Big) 
\ll \Big(\frac{\pi k}{50}\Big)^{k-1} \frac{1}{(k-2)!} e^{\pi k/50} \ll e^{-k}, 
$$
for large $k$, as desired.  

\enddemo 

Since $n_1 \cdots n_r \le x^r = \sqrt{k}$ we see by Lemma 2.1 that
$$
S_2= \sum\Sb n_1, \ldots, n_r \le x \endSb \frac{b_1(n_1,\ldots,n_r)}{\sqrt{n_1\cdots n_r}}
+O\Big( e^{-k} \sum_{n_1,\ldots n_r\le x} 
\frac{\tau(n_1)\cdots \tau(n_r)\tau(n_1\cdots n_r)}{\sqrt{n_1 \cdots n_r}}\Big).
$$
The error term is easily seen to be $\ll e^{-k} x^k = k^{\frac 12} e^{-k}$, 
a negligible amount.  As for the main term we see easily that 
$$
\sum_{n\le x} \frac{B_r(n)}{\sqrt{n}} \le \sum\Sb n_1, \ldots, n_r \le x \endSb 
\frac{b_1(n_1,\ldots,n_r)}{\sqrt{n_1\cdots n_r}} 
\le \sum_{n\le x^r} \frac{B_r(n)}{\sqrt{n}}.  
$$
Recall that $B_r(n)$ is a multiplicative function with $B_r(p)=0$, 
$B_r(p^2)=r(r-1)/2$ and $B_r(p^a)$ grows only polynomially in $a$.  
Thus the generating function $\sum_{n=1}^{\infty} B_r(n)n^{-s}$ 
can be compared with $\zeta(2s)^{r(r-1)/2}$, the quotient being 
a Dirichlet series absolutely convergent in Re$(s)>\frac{1}{4}$.  
A standard argument 
(see Theorem 2 of [9]) therefore gives that 
$$
\sum_{n\le z} \frac{B_r(n)}{\sqrt{n}} \sim C_r (\log z)^{r(r-1)/2}, 
$$
for a positive constant $C_r$.   It follows that 
$$
S_2 \asymp (\log x)^{r(r-1)/2} \asymp (\log k)^{r(r-1)/2}.
$$

We now turn to $S_1$.  To evaluate $S_1$ we need an `approximate 
functional equation' for $L(\frac 12,f)$.  

\proclaim{Lemma 2.2} Define for any positive number $\xi$ the 
weight 
$$
W_k(\xi):= \frac{1}{2\pi i} \int_{(c)} \frac{(2\pi)^{-\frac 12 -s} 
\Gamma(s+\frac{k}{2})}{(2\pi)^{-\frac 12} \Gamma(\frac k2)} \xi^{-s} \frac{ds}{s},
$$
where the integral is over a vertical line $c-i\infty$ to $c+i\infty$ with 
$c>0$.  Then, for $k\equiv 0\pmod 4$, 
$$
L(\tfrac 12,f) = 2 \sum_{n=1}^{\infty} \frac{\lam_f(n)}{\sqrt{n}} W_k(n) 
= 2 \sum_{n\le k} \frac{\lam_f(n)}{\sqrt{n}} W_k(n) + O(e^{-k}).
$$
Further the weight $W_k(\xi)$ satisfies $|W_k(\xi)|\ll k\pi^{-k}/\xi$ 
for $\xi >k$, $W_k(\xi)= 1+ O(e^{-k})$ for $\xi<k/100$, and 
$W_k(\xi) \ll 1$ for $k/100 \le \xi \le k$.   
\endproclaim 
\demo{Proof}  The argument is standard.  For $1\le c>\tfrac 12$ we consider 
$$
I =\frac{1}{2\pi i} \int_{(c)} \frac{(2\pi)^{-\frac 12-s} \Gamma(s+\frac k2)}{(2\pi)^{-\frac 12} 
\Gamma(\frac k2)} L(\tfrac 12+s,f)\frac{ds}{s}.
$$
Expanding out $L(\tfrac 12+s,f)$ and integrating term by term we see 
that 
$$I= \sum_{n=1}^{\infty} \lam_f(n)n^{-\frac 12} W_k(n).
$$  
On the other hand moving the line of integration to the line Re$(s)=-c$ we 
see that 
$$
I= L(\tfrac 12,f) + \frac{1}{2\pi i} \int_{(-c)} \frac{\Lam(\tfrac 12+s,f)}{(2\pi)^{-\frac 12} \Gamma(
\frac k2)} \frac{ds}{s},
$$
and using the functional equation $\Lam(\tfrac 12+s,f)=\Lam(\tfrac 12-s,f)$, and 
replacing $s$ by $-s$ in the integral above, we see that $I=L(\tfrac 12,f)-I$.  Thus 
$$
L(\tfrac 12,f) =2I  = 2\sum_{n=1}^{\infty} \frac{\lam_f(n)}{\sqrt{n}} W_k(n). 
$$

Regarding the weight $W_k(\xi)$ note that by considering the integral 
for some large positive integer $c$ we get that
$$
\align
|W_k(\xi)| &\le \frac{1}{2\pi} \int_{(c)} (2\pi \xi)^{-c} \frac{|\Gamma(s+\frac k2 +1)|}{\Gamma(\frac k2)} 
\frac{|ds|}{|s(s+k/2)|} 
\\
&\le (2\pi \xi)^{-c} \frac{\Gamma(c+1+\frac{k}{2})}{\Gamma(\frac k2)} 
\le (2\pi \xi)^{-c} (k+c)^c.\\
\endalign
$$
Taking $c=k$ we obtain that $|W_k(\xi)| \le (k/(\pi \xi))^k$ 
so that if $\xi \ge k$ then $|W_k(\xi)| \le (k/\xi) \pi^{-(k-1)}$.  This 
proves the first bound for $W_k(\xi)$ claimed in the Lemma, and also shows that 
$$
\Big|\sum_{n> k } \frac{\lam_f(n)}{\sqrt{n}} W_k(n) \Big| 
\ll k \pi^{-k} \sum_{n>k} \frac{|\lam_f(n)|}{n^{\frac 32}} 
\ll e^{-k}.
$$
The other claims on $W_k(\xi)$ are proved similarly; for the range $\xi <k/100$ 
we move the line of integration to $c=-\frac k2 +1$, for the last range 
$k/100 \le \xi \le k$ just take the integral to be on the line $c=1$.  
\enddemo 

Returning to $S_1$ note that 
$$
A(f)^{r-1} = \sum\Sb n_1, \ldots, n_{r-1}\le x \endSb 
\frac{1}{\sqrt{n_1 \cdots n_{r-1}}} \sum_{t|n_1\cdots n_{r-1}} b_t(n_1,\ldots,n_{r-1}) 
\lam_f(t).
$$
Since $A(f)^{r-1}$ is trivially seen to be $\ll x^{r-1} <\sqrt{k}$, we see 
by Lemma 2.2, that 
$$
S_1 =2 \sum_{n\le k} \frac{1}{\sqrt{n}} W_k(n) \sum_{n_1, \ldots, n_{r-1}\le x} 
\sum_{t|n_1 \cdots n_{r-1}} \frac{b_{t}(n_1,\ldots,n_{r-1})}{\sqrt{n_1\cdots n_{r-1}} }
\sumh_{f\in {\Cal H}_k} \lam_f(t)\lam_f(n) + O(\sqrt{k}e^{-k}). 
$$
Now we appeal to Lemma 2.1.  The error term that arises is 
trivially bounded by $\ll ke^{-k}$ which is negligible.  In the 
main term $\delta(n,t)$, since $t\le x^{r-1}<\sqrt{k}$ 
we may replace $W_k(n)$ by $1+O(e^{-k})$.  It follows that 
$$
S_1 = 2 \sum\Sb n_1,\ldots, n_{r-1} \le x \endSb \sum_{t|n_1\cdots n_{r-1}} \frac{b_t(n_1,\ldots,n_{r-1})}
{\sqrt{n_1\cdots n_{r-1}}} \frac{1}{\sqrt{t}}+ O(ke^{-k}). 
$$
Now observe that $b_1(n_1,\ldots,n_{r-1},t)=b_t(n_1,\ldots,n_{r-1})$ 
if $t$ divides $n_1\cdots n_{r-1}$, and otherwise $b_1(n_1,\ldots,n_{r-1},t)$ is 
zero.  Therefore, writing $n_r$ for $t$, we obtain that 
$$
S_1 = 2\sum\Sb n_1 ,\ldots, n_{r-1} \le x \endSb \sum_{n_r \le \sqrt{k}} 
\frac{b_1(n_1,\ldots,n_r)}{\sqrt{n_1\cdots n_r}} +O(ke^{-k}).
$$
Using $b_1\geq 0$ we see that $S_1 \ge 2 S_2 +O(ke^{-k})$, and moreover, 
arguing as in the case of  $S_2$ we may see that 
$$
S_1\asymp (\log k)^{r(r-1)/2}.
$$
Theorem 1 follows.

\head 3. Proof of Theorem 2 \endhead

\noindent For simplicity, we will restrict ourselves to fundamental discriminants 
of the form $8d$ where $d$ is a positive, odd square-free number 
with $X/16 < d\le X/8$.  %Then the sum over these is 
%$$ 
%\sumf_{d\sim X} c_d = \sum_{X/16<d<X/8} \mu(2d)^2 c_d .
%$$
Let $k$ be a given even number, and set 
$x= X^{\frac{1}{10k}}$.  Define 
$$
A(8d):= \sum_{n \le x} \frac{\chi_{8d}(n)}{\sqrt{n}}, 
$$
and let 
$$
S_1 := \sum\Sb X/16 < d\le X/8 \endSb \mu^2 (2d) L(\tfrac 12, \chi_{8d})A({8d})^{k-1},
\ \  \text{and} \ \ 
S_2:=\sum\Sb X/16 < d\le X/8 \endSb \mu^2(2d) A({8d})^{k}.
$$
An application of H{\" o}lder's inequality 
gives that 
$$
\sumf_{|d|\le X} L(\tfrac 12,\chi_{8d})^k  \ge 
 \sum_{X/16<d\le X/8} \mu^2(2d) L(\tfrac 12,\chi_{8d})^{k} 
\ge \frac{S_1^{k}}{S_2^{k-1}},
$$
so that to prove Theorem 2 we need only give satisfactory estimates
for 
$S_1$ and $S_2$.  

\def\L{\fracwithdelims()}
We start with $S_2$.  Expanding our $A(8d)^k$ we see that 
$$
S_2 = \sum_{n_1,\ldots,n_k \le x} \frac{1}{\sqrt{n_1\cdots n_k}} 
\sum\Sb X/16 < d\le X/8 \endSb \mu^2(2d) \L{8d}{n_1 \cdots n_k }. \tag{3.1} 
$$

\proclaim{Lemma 3.1}   Let $n$ be an odd integer, and let $z\ge 3$ be a real 
number.  If $n$ is not a perfect square then
$$
\sum\Sb d\le z \endSb \mu^2(2d) \L{8d}{n} \ll z^{\frac 12} n^{\frac 14} \log (2n),  
$$
while if $n$ is a perfect square then 
$$
\sum\Sb d\le z\endSb \mu^2(2d) \L{8d}{n} =  \frac{z}{\zeta(2)} \prod_{p|2n} \L{p}{p+1}+ O(z^{\frac 12+\epsilon}n^{\epsilon}).
$$
\endproclaim 
\demo{Proof}  Note that $\sum_{\alpha^2|d} \mu(\alpha)=1$ if $d$ is square-free 
and $0$ otherwise.  Therefore, 
$$
\sum_{n\le z} \mu^2(2d) \L{8d}{n} 
= \sum\Sb \alpha \le \sqrt{z} \\ \alpha \text{ odd} \endSb \L{8\alpha^2}{n}
\sum\Sb d\le z/\alpha^2 \\ d\text{ odd} \endSb  \L{d}{n}.
$$
If $n$ is not a square then the inner sum over $d$ is a character sum 
to a non-principal character of modulus $2n$ (we take $2n$ to account for 
$d$ being odd), and the P{\' o}lya-Vinogradov inequality (see [2]) gives 
that the sum over $d$ is $\ll \sqrt{n} \log (2n)$.  Further, the 
sum over $d$ is trivially $\ll z/\alpha^2$.  Thus, if $n$ is not a square, we 
get that 
$$
\sum_{n\le z} \mu^2(2d) \L{8d}{n} \ll \sum\Sb \alpha \le \sqrt{z} \endSb 
\min \Big( \sqrt{n} \log (2n), \frac{z}{\alpha^2}\Big) 
\ll z^{\frac 12} n^{\frac 14} \log (2n),
$$
upon using the first bound for $\alpha \le z^{\frac 12} n^{-\frac 14}$ 
and the second bound for larger $\alpha$.  

If $n$ is a perfect square, then $\L{8d}{n}=1$ if $d$ is coprime to $n$, 
and is $0$ otherwise.  Thus 
$$
\sum_{d\le z} \mu(2d)^2 \L{8d}{n} = 
\sum\Sb d\le z \\ (d,2n)=1 \endSb \mu^2(d) = \frac{z}{\zeta(2)} \prod_{p|2n} \L{p}{p+1}  +O(z^{\frac 12+\epsilon}n^{\epsilon}),
$$ 
by a standard argument.

\enddemo

Using Lemma 3.1 in (3.1) we obtain that 
$$ 
S_2 = \frac{X}{16 \zeta(2)} \sum\Sb n_1, \ldots , n_k \le x\\ 
n_1 \cdots n_k = \text{ odd square} \endSb 
\frac{1}{\sqrt{n_1 \cdots n_k}} \prod_{p|2 n_1\cdots n_k} \L{p}{p+1} 
+ O\Big( X^{\frac 12+\epsilon} x^{k(\frac 34+\epsilon)}\Big).
$$ 
Since $x=X^{\frac 1{10k}}$ the error term above is $\ll X^{\frac 35}$.  Writing 
$n_1 \cdots n_k =m^2$ we see that 
$$
\align
\sum\Sb m^2 \le x \\ m\text{ odd} \endSb 
\frac{d_k(m^2)}{m} \prod_{p|2m} \L{p}{p+1} 
&\le \sum\Sb n_1, \ldots , n_k \le x\\ 
n_1 \cdots n_k = \text{ odd square} \endSb 
\frac{1}{\sqrt{n_1 \cdots n_k}} \prod_{p|2 n_1\cdots n_k} \L{p}{p+1} 
\\
&\le \sum\Sb m^2 \le x^k \\ m\text{ odd} \endSb \frac{d_k(m^2)}{m} \prod_{p|2m} \L{p}{p+1}. \\
\endalign
$$
A standard argument (see Theorem 2 of [9]) shows that 
$$
\sum\Sb m\le z \\ m\text{ odd} \endSb  \frac{d_k(m^2)}{m} \prod_{p|2m} \L{p}{p+1} 
\sim C(k) (\log z)^{k(k+1)/2},
$$
for a positive constant $C(k)$.  We conclude that 
$$
S_2 \asymp X (\log X)^{k(k+1)/2}. \tag{3.2}  
$$ 

It remains to evaluate $S_1$.  As before, we need an `approximate 
functional equation' for $L(\tfrac 12,\chi_{8d})$.  

\proclaim{Lemma 3.2}   For a 
positive number $\xi$ define the weight 
$$
W(\xi) = \frac{1}{2\pi i} \int_{(c)} \frac{\Gamma(\frac s2 +\frac 14)}{\Gamma(\frac 14)} 
\xi^{-s} \frac{ds}{s},
$$
where the integral is over a vertical line $c-i\infty$ to $c+i\infty$ with $c>0$.  Then, 
for any odd, positive, square-free number $d$ we have 
$$
L(\tfrac 12,\chi_{8d}) 
= 2\sum_{n=1}^{\infty} \frac{\chi_{8d}(n)}{\sqrt{n}} W\Big(\frac{n\sqrt{\pi}}{\sqrt{8d}}\Big).
$$
The weight $W(\xi)$ is smooth and satisfies $W(\xi) = 1+ O(\xi^{\frac 12-\epsilon})$ 
for $\xi$ small, and for large $\xi$ satisfies $W(\xi) \ll e^{-\xi}$.  Moreover 
the derivative $W^{\prime}(\xi)$ satisfies $W^{\prime}(\xi) \ll \xi^{\frac 12-\epsilon} e^{-\xi}$.  
\endproclaim 
\demo{Proof}  This is given in Lemmas 2.1 and 2.2 of [10], but 
for completeness we give a sketch.  For some $1\ge c > \frac 12$, we consider 
$$
\frac{1}{2\pi i} \int_{c)} \frac{(8d/\pi)^{\frac s2 + \frac 14} \Gamma(\frac s2+\frac 14)} {(8d/\pi)^{\frac 14} 
\Gamma(\frac 14)} L(\tfrac 12+s,\chi_{8d}) \frac{ds}{s}, 
$$
and argue exactly as in Lemma 2.2.   This gives the desired formula for 
$L(\tfrac 12,\chi_{8d})$.  The results on the weight $W(\xi)$ follow upon moving 
the line of integration to Re$(s)=-\frac 12+\epsilon$ when $\xi$ is small, and taking $c$ to  
be an appropriately large positive number if $\xi$ is large.  
\enddemo

By Lemma 3.2 we see that 
$$
S_1= 2 \sum_{n=1}^{\infty} \frac{1}{\sqrt{n}} \sum_{n_1,\ldots, n_{k-1}\le x} 
\frac{1}{\sqrt{n_1\cdots n_{k-1}}} 
\sum\Sb X/16 < d\le X/8\endSb \mu^2(2d) \L{8d}{n n_1\cdots n_{k-1}} 
W\L{n\sqrt{\pi}}{\sqrt{8d}}.\tag{3.3} 
$$
If $n n_1 \cdots n_{k-1}$ is not a square, then using Lemma 3.1 and 
partial summation we may see that 
$$
 \sum\Sb X/16 < d\le X/8\endSb \mu^2(2d) \L{8d}{n n_1\cdots n_{k-1}} 
W\L{n\sqrt{\pi}}{\sqrt{8d}} 
\ll X^{\frac 12} (nn_1\cdots n_{k-1})^{\frac 14 +\epsilon}  e^{-n/\sqrt{X}}.
$$
If $n n_1\cdots n_{k-1}$ is an odd square then Lemma 3.1 and partial 
summation gives that the sum over $d$ in (3.3) is 
$$
\frac{X}{16 \zeta(2)} \prod_{p| 2 nn_1\cdots n_{k-1}} 
\L{p}{p+1} \int_1^2 W\L{n\sqrt{2\pi}}{\sqrt{Xt}} dt + O(X^{\frac 12+\epsilon} e^{-n/\sqrt{X}}).
$$
We use these two observations in (3.3).  Note that the error terms 
contribute to (3.3) an amount $\ll X^{\frac 12+\epsilon} x^{(k-1)(\frac 34+\epsilon)} X^{\frac 38+\epsilon} 
\ll X^{\frac {39}{40}}$.  It remains to estimate the main term contribution to (3.3).  
To analyze these terms let us write $n_1\cdots n_{k-1}$ as $r s^2$ where $r$ and $s$ are 
odd and $r$ is square-free.  Then $n$ must be of the form $r\ell^2$ where $\ell$ is odd.  
With this notation the main term contribution to (3.3) is 
$$
\frac{X}{8\zeta(2)} \sum\Sb  rs^2 =n_1 \cdots n_{k-1} \\ n_1, \ldots, n_{k-1} \le x \endSb 
\frac{1}{rs} \sum\Sb \ell \text{ odd } \endSb 
\frac{1}{\ell} \int_1^2 \prod_{p|2rs\ell} \L{p}{p+1} W\Big( \frac{r\ell^2 \sqrt{2\pi}}{\sqrt{Xt}}\Big) dt.
$$
Note that $r \le x^{k-1} < X^{\frac 1{10}}$, and an easy calculation gives  
that  the sum over $\ell$ above is 
$$
= \prod_{p|2rs} \L{p}{p+1} \prod_{p\nmid 2rs} 
\Big(1-\frac{1}{p(p+1)} \Big) \frac{1}{4} 
\log \frac{\sqrt{X}}{r} + O(1).
$$
It follows that the main term contribution to (3.3) is 
$$
\align
&\gg X\log X \sum\Sb rs^2= n_1\cdots n_{k-1} \\ n_1, \ldots, n_{k-1} \le x 
\endSb \frac{1}{rs} \prod_{p|2rs} \L{p}{p+1} 
\\
&\gg X\log X \sum\Sb r \text{ odd and square-free}\\ s\text{ odd} \\ rs^2 \le x \endSb \frac{d_{k-1}(rs^2)}{rs} 
\prod_{p|2rs} \L{p}{p+1}\gg X (\log X) (\log x)^{k-1+k(k-1)/2},
\\
\endalign
$$
where the last bound follows by invoking Theorem 2 of [9].   
We conclude that 
$$
S_1 \gg X(\log X)^{k(k+1)/2},
$$
which when combined with (3.2) proves Theorem 2.

 \enddocument